\documentclass[imslayout,preprint,11pt]{imsart}
\RequirePackage[OT1]{fontenc}
\RequirePackage{amsthm,amsmath}
\RequirePackage{natbib}
\startlocaldefs

\numberwithin{equation}{section}
\theoremstyle{plain}
\newtheorem{thm}{Theorem}[section]
\newtheorem{lemma}{Lemma}[section]
\newtheorem{corollary}{Corollary}[section]

\endlocaldefs

\begin{document}

\begin{frontmatter}
\title{Some notes on improving upon the James-Stein estimator 
}
\runtitle{Improving on the JSE}

\begin{aug}
\author{\fnms{Yuzo} \snm{Maruyama}
\ead[label=e1]
{maruyama@csis.u-tokyo.ac.jp}}

\address{
Center for Spatial Information Science, \\
The University of Tokyo\\
5-1-5 Kashiwanoha, Kashiwa-shi, Chiba
277-8568, Japan\\
\printead{e1}}



\runauthor{Y. Maruyama}

\affiliation{The University of Tokyo}

\end{aug}

\begin{abstract}
We consider estimation of a multivariate normal mean vector under sum
of squared error loss. We propose a new class of 
smooth estimators parameterized by $\alpha$
dominating the James-Stein estimator.
The estimator for $\alpha=1$ corresponds to the generalized Bayes
estimator with respect to the harmonic prior.
When $ \alpha $ goes to infinity, the estimator converges to 
the James-Stein positive-part estimator.
Thus the class of our estimators is a bridge between the admissible estimator
($\alpha=1$) and the inadmissible estimator ($\alpha=\infty$).
Although the estimators have quasi-admissibility which is 
a weaker optimality than admissibility, 
the problem of determining whether or not the estimator for $\alpha>1$
admissible is still open.
\end{abstract}


\end{frontmatter}

\section{Introduction}
Let $ X $ be a random variable having $ p $-variate normal distribution
$ N_{p} (\theta,I_{p}) $. Then we consider the problem of estimating  
the mean vector $ \theta $ by $ \delta(X) $ relative to  quadratic
loss.
Therefore every estimator is evaluated based on the risk function
\begin{equation*}
 R(\theta,\delta)=E_{\theta}\left[ \| \delta(X) - \theta \|
 ^{2} \right]=\int_{R^{p}}\frac{\| \delta(x) - \theta \|^{2}}{(2\pi)^{p/2}}
\exp\left(-\frac{\| x - \theta \|^{2}}{2}\right)dx.  
\end{equation*}
The usual estimator $ X $, with the constant risk $ p $, is minimax for
any dimension $p$.
It is also admissible when $p=1$ and $2$, as shown in \cite{Blyth1951}
 and \cite{Stein1956}, respectively. 
\cite{Stein1956} showed, however, that 
when $p \geq 3$, 
there exists an estimator dominating $ X $ among a class of
equivariant estimators 
relative to the orthogonal transformation group which have the form  
\begin{align} \label{estiortho}
 \delta_{\phi}(X)= ( 1- \phi(\| X \|^{2})/\| X 
\|^{2} )X.
\end{align}
\cite{James-Stein1961} succeeded in giving an explicit form of an
 estimator improving on $ X $ as 
\begin{equation*}
\delta_{JS}(X)=(1-(p-2)/\| X \|^{2})X,
\end{equation*}  
which is called the James-Stein estimator.
More generally, a  large class of better estimators than $X$ 
has been proposed in the literature.
The strongest tool for this is \cite{Stein1974}'s identity as follows.
\begin{lemma}[\cite{Stein1974}] \label{SI}
If $ Y \sim  N(\mu, 1) $ and $ h(y) $ is any differentiable
 function such that $E[|h'(Y)|] < \infty $, then
\begin{equation} 
 E[ h(Y)(Y-\mu)]= E[ h'(Y)].  \label{stein}
\end{equation}
\end{lemma}
By using the identity \eqref{stein}, 
 the risk function of the estimator of
the form $ \delta_g(X)=X+g(X)=(X_1+g_1(X), \dots , X_p+g_p(X))' $ is written as
\begin{eqnarray*}
 R(\theta,\delta_g)&=& E_{\theta}\left[
\|\delta_g(X)-\theta\|^{2} \right]\nonumber \\
&=& E_{\theta}\left[ \|X-\theta\|^{2}\right]
+E_{\theta}\left[\|g(X)\|^{2}\right]
+2\sum_{i=1}^{p}
E_{\theta}\left[(X-\theta)'g(X)\right] \nonumber \\ 
 &=& p+E_{\theta}\left[\|g(X)\|^{2}\right]+2\sum_{i=1}^{p}
E_{\theta}\left[ \frac{\partial}{\partial x_{i}}g_{i}(X)\right], 
\label{fuhenrisk} 
\end{eqnarray*}
where  $ 
g_{i}(x)  $ is assumed to be differentiable and  
$ E|(\partial/\partial
x_{i})g_{i}(X)| < \infty $  
for $ i=1,\ldots, p. $
Since a statistic $\hat{R}(\delta(X))$ given by
\begin{align}\label{Unbias}
 p+\|g(X)\|^{2}+2\sum_{i=1}^{p}
\frac{\partial}{\partial x_{i}}g_{i}(X)
\end{align}
does not depend on the unknown parameter $\theta$ and satisfies
$E(\hat{R}(\delta(X)))=R(\theta,\delta)$,
it is called  the unbiased estimator of risk.
Clearly if 
\begin{equation} \label{parcon}
 \|g(x)\|^{2}+2\sum_{i=1}^{p}
\frac{\partial}{\partial x_{i}}g_{i}(x) \leq 0,
\end{equation}
for any $x$, then $ \delta_g(X) $ 
 dominates $X$.

To make the structure of  estimators improving on $X$
more comprehensible, 
we consider a class of orthogonally equivariant estimators of a form 
given in \eqref{estiortho}.
Assigning   $  g(x)=
-\phi(\|x\|^{2})x/\|x\|^{2} $ in \eqref{parcon},
\eqref{parcon} becomes
\begin{equation} 
 \phi(w)\left( 2(p-2)-\phi(w)\right)/w+4\phi'(w) \geq 0, \label{diffcon}
\end{equation}
for $w=\|x\|^2$. 
The inequality \eqref{diffcon} is, for example, satisfied if
$\phi(w)$ is monotone nondecreasing and 
within $[0,2(p-2)]$ for any $w\geq 0$.

Since  $ \phi(w)=c $ for
$0<c<2(p-2)$
 satisfies the inequality \eqref{diffcon},
\begin{align} \label{jsc}
 (1-c/\|X\|^2)X
\end{align}
with $0<c<2(p-2)$ dominates $X$.
The estimator $ \delta_{JS} $, 
\eqref{jsc} with $c=p-2$,
is the best estimator
among a class of estimators \eqref{jsc} because
the risk function of \eqref{jsc} is given by
\begin{align*}
p+c(c-2(p-2))E[\|X\|^{-2}] 
\end{align*}
and hence minimized by $c=p-2$.


It is however noted that
when $ \|x \|^{2} < p-2 $, the James-Stein estimator yields an
over-shrinkage and changes the sign of each component of $ X $. 
The James-Stein positive-part estimator
\begin{equation*}
\delta_{JS}^{+}(X)=\max(0,1-(p-2)/\| X \|^{2})X,
\end{equation*}
eliminates this drawback and 
dominates the James-Stein
estimator.
We notice here that 
the technique for proving the
inadmissibility of $\delta_{JS}$ is not from the Stein identity or
the unbiased estimator
of risk given in \eqref{Unbias}.
The risk difference between $\delta_{JS}$ and $\delta_{\phi}$ is 
given by
\begin{align*}
 R(\delta_{JS},\theta)-R(\delta_{\phi},\theta)
=E\left[-\frac{(\phi(\|X\|^2)-p+2)^2}{\|X\|^2}+4\phi'(\|X\|^2)
\right],
\end{align*}
but $\phi_{JS}^+(w) =\min(w,p-2)$ which makes the James-Stein
positive-part estimator does not satisfy the
inequality
\begin{align} \label{js-ineq}
 -(\phi(w)-p+2)^2+4w\phi'(w) \ge 0
\end{align}
for any $w \ge 0$.
In Section 2, we will show that there is no $\phi $ which satisfies
\eqref{js-ineq} for any $w \ge 0$, that is,
the Stein identity itself is not useful 
for finding estimators dominating $\delta_{JS}$.
We call such optimality for the James-Stein estimator
quasi-admissibility.
We will explain the concept and give a sufficient condition for
quasi-admissibility in Section 2.

In spite of such difficulty, some 
estimators which dominate the James-Stein estimator have been given
by several authors.
\cite{Li-Kuo1982} and \cite{Guo-Pal1992} considered 
the class of estimators of forms
$ \delta_{LK}(X)=(1-\phi_{LK}(\|X\|^2)/\|X\|^2)X $ where
\begin{equation*}
 \phi_{LK}(w)=p-2-\sum_{i=1}^{n}a_iw^{-b_i}, 
\end{equation*}
where $ a_i \geq 0 $ for any $ i$ and $ 0 < b_1 < b_2 <
  \dots < b_n$.
For example when $n=1$, they both showed that, 
$ \delta_{LK}(X)$ for $ 0 < b_1 < 4^{-1}(p-2) $ and
  $a_1=2b_12^{b_1}\Gamma(p/2-b_1-1)/\Gamma
(p/2-2b_1-1)$ is superior to the James-Stein estimator.  
\cite{Kuriki-Takemura2000} gave two estimators
which shrink toward the ball with center $0$, 
 $ \delta_{KT}^i(X)=(1-\phi_{KT}^{i}(\|X\|^2)/\|X\|^2)X $ for $i=1,2$
where
\begin{equation*}
 \phi_{KT}^1(w)=\begin{cases}
		 0 & w \leq r^2 \\
		 p-2 -\sum_{i=1}^{p-2}(r/w^{1/2})^i & w > r^2,
		\end{cases}
\end{equation*}
\begin{equation*}
 \phi_{KT}^2(w)=\begin{cases}
	      0 & w \leq \{ (p-1)/(p-2)\}^2 r^2 \\
	      p-2-r/(w^{1/2}-r) & w > \{ (p-1)/(p-2)\}^2 r^2.
	     \end{cases}
\end{equation*}
They showed that when $r$ is sufficiently small,
these two estimators dominate the James-Stein estimator.
However, these estimators are not appealing 
since they are inadmissible.
In our setting, the estimation of a multivariate normal mean,
\cite{Brown1971} showed that any admissible estimator should be generalized
Bayes. Since the shrinkage factor $(1-\phi_{LK}(w)/w)$ becomes negative for some
$w$ and $ \phi_{KT}^i$ for $i=1,2$ fail to be analytic,
neither $ \delta_{LK}$ nor $ \delta_{KT}^{i}$ can be generalized Bayes
or admissible.

In general, when 
we propose an estimator ($\delta_*$, say) dominating a certain inadmissible
estimator, it is extremely important to find it among admissible
estimators. If not, a more difficult problem (finding an estimator 
improving on $\delta_*$)
just occurs.
To the best of our knowledge, 
the sole admissible estimator dominating the James-Stein
estimator is \cite{Kubokawa1991}'s estimator
$ \delta_{K}(X)=(1-\phi_{K}(\|X\|^2)/\|X\|^2)X$ where
\begin{equation}\label{kubo-phi}
  \phi_{K}(w)=p-2-2\frac{\exp(-w/2)}
{\int_{0}^{1}
\lambda^{p/2-2}\exp(-w\lambda/2) d \lambda}.
\end{equation} 
The estimator is  generalized Bayes
with respect to the harmonic prior density 
$ \| \theta \|^{2-p} $ which was originally suggested by \cite{Stein1974}.
The only fault is, however, that
it does not improve upon $ \delta_{JS}(X) $ at $ \| \theta \|=0 $.
(See Section 3 for the detail.)
Shrinkage estimators like \eqref{estiortho}
make use of the vague prior information
that $ \| \theta \| $ is close to $ 0 $.
It goes without saying that we would like to get the 
significant improvement 
of risk when the prior information is accurate. 
Though $ \delta_{K}(X) $ is an admissible generalized Bayes estimator 
and thus smooth, it has no improvement on $ \delta_{JS}(X) $ at the
origin $\| \theta \|=0$. 
On the other hand, $ \delta_{JS}^+(X) $ 
improves on $ \delta_{JS}(X) $ 
significantly at the origin,
but it is not analytic and is thus inadmissible by \cite{Brown1971}'s
complete class theorem.
Therefore a more challenging open problem is to find admissible
estimators dominating the James-Stein estimator especially at 
$ \|\theta \|=0 $. 
In this paper, we 
will consider a class of estimators
\begin{equation} 
 \delta_{\alpha}(X)=
\frac{\int_{0}^{1}(1-\lambda) \lambda^
{ \alpha(p/2-1)-1}\exp(-\|X\|^2 \alpha \lambda/2
)d\lambda}{\int_{0}^{1}\lambda^
{ \alpha(p/2-1)-1}\exp(-\|X\|^2 \alpha \lambda/2)d\lambda}
X
\end{equation}
for $\alpha \ge 1$.
In Section 3, we show that $ \delta_\alpha $ with $ \alpha \geq 1$
improves on the James-Stein estimator and
that it has strictly risk improvement at $\|\theta\|=0$.
Furthermore we see that
$ \delta_\alpha $ approaches $\delta_{JS}^+ $ as $\alpha $ goes to $\infty$.
Since $\delta_\alpha$ with $ \alpha=1$ corresponds to $\delta_K$,
the class of $ \delta_\alpha $ with $ \alpha \ge 1$ is a bridge between
$\delta_K$ which is admissible and $ \delta_{JS}^+$ which is
inadmissible. Although we show that $ \delta_\alpha $ with $ \alpha > 1$ 
is quasi-admissible, which is introduced in Section 2,
we have no idea on its admissibility at this stage.

\section{Quasi-admissibility}
In this section, we introduce the concept of quasi-admissibility
and give a sufficient condition for quasi-admissibility.
We deal with a reasonable class of estimators which have the form
\begin{align}
\delta_m &= X+\nabla \log m(\|X\|^{2})  \label{pseudo} \\
&=\{1+2m'(\|X\|^2)/m(\|X\|^2)\}X \notag
\end{align}
where $ \nabla $ is a differential operator
$
( \partial/\partial x_1, \dots, \partial/\partial x_p)'
$
and $m$ is a positive function.
If $m(w)=w^{-c} $ for $c>0$, \eqref{pseudo}
becomes the James-Stein type estimator \eqref{jsc}.
Any (generalized) Bayes estimator
with respect to spherical symmetric measure $\pi$
should also have the form \eqref{pseudo} because
it is written as
\begin{align*}
& \frac{\int_{R^p} \theta \exp(-\|x - \theta\|^2/2)\pi(d\theta)}
{\int_{R^p} \exp(-\|x - \theta\|^2/2)\pi(d\theta)} \\
& = x + \frac{\int_{R^p} (\theta -x) \exp(-\|x - \theta\|^2/2)\pi(d\theta)}
{\int_{R^p} \exp(-\|x - \theta\|^2/2)\pi(d\theta)}\\
& = x +  \frac{\nabla \int_{R^p} \exp(-\|x - \theta\|^2/2)\pi(d\theta)}
{\int_{R^p} \exp(-\|x - \theta\|^2/2)\pi(d\theta)} \\
& = x +  \nabla \log \int_{R^p} \exp(-\|x -
 \theta\|^2/2)\pi(d\theta) \\
& = x +  \nabla \log m_\pi(\|x\|^2) 
\end{align*}
where $ m_\pi(\|x\|^2)=\int_{R^p} \exp(-\|x -
\theta\|^2/2)\pi(d\theta) $.
\cite{Bock1988} and \cite{Brown-Zhang2006}
called an estimator of the form \eqref{pseudo} pseudo-Bayes
and quasi-Bayes, respectively.
If, for given $m$, there exists a nonnegative measure $ \nu $ which satisfies 
\begin{equation} \label{sokudo}
m(\|x\|^{2})=\int_{\mathit{R^{p}}} 
\exp (-\| \theta -x\|^{2}/2) \nu(d\theta), 
\end{equation}
the estimator of the form \eqref{pseudo} is truly generalized
Bayes. 
However it is often difficult to 
determine whether or not 
$m$ has such an exact integral form.

Suppose that 
$ \delta_{m,k}(X)=\delta_m(X)+k(\|X\|^2)X$ is a competitor of $ \delta_m$.
%
Then substituting $ g(x)=\nabla \log m(\| x \|^{2}) =2 x
m'(\|x\|^2)/m(\|x\|^2)$ and  $ g(x)=
2x m'(\|x\|^2)/m(\|x\|^2)+ k(\|x\|^2)x $ in 
\eqref{Unbias} respectively, we have the unbiased estimators of risk as 
\begin{equation} 
\hat{R}(\delta_{m})= p
-4w\left(\frac{m'(w)}{m(w)}\right)^{2}+4p\frac{m'(w)}{m(w)}
+8w\frac{m''(w)}{m(w)}
\end{equation}
for $ w=\|x\|^{2}$ 
and $ \hat{R}(\delta_{m,k})= \hat{R}(\delta_{m})+ \Delta(m,mk)$ where
\begin{equation} \label{Delta}
\Delta(m,mk)=
4w\frac{m'(w)}{m(w)}k(w)+2pk(w)+4wk'(w)+wk^2(w).
\end{equation}
If there exists $ k $ such that $ \Delta(m,mk) \leq 0$ for any
 $w \geq 0$ with strict inequality for some $ w $
, that implies that 
$ \delta_m$ is inadmissible, that is,
$ R(\theta,\delta_{mk}) \leq R(\theta,\delta_m)$
for all $\theta$  with strict inequality for some $ \theta $.
If there does not exist such $k$,
$ \delta_m$ is said to be quasi-admissible.
Hence quasi-admissibility is a weaker optimality than admissibility.
Now we state a sufficient condition for quasi-admissibility. The idea is 
originally from \cite{Brown1988}, but the paper is not so accessible.
See also \cite{Rukhin1995}, where quasi-admissibility is called
 permissibility. 
%
%

%
\begin{thm} \label{perm}
The estimator of the form \eqref{pseudo} is quasi-admissible if 
\begin{align*}
\int_{0}^{1}w^{-p/2}m(w)^{-1}dw = \infty \ \mbox{ as well as } \
 \int_{1}^{\infty}w^{-p/2}m(w)^{-1}dw = \infty.
\end{align*}
\end{thm}
\begin{proof}
We have only to show that $\hat{R}(\delta_{m,k}) 
\leq  \hat{R}(\delta_{m})  $, that is, $\Delta(m,mk)\leq 0$
 implies $ k(w) \equiv 0 $.
Let $M(w)=w^{p/2}m(w)$ and $h(w)=M(w)k(w)$.
Then we have
\begin{align*}
 \Delta(m,mk)= 4w\frac{h'(w)}{M(w)}+w\frac{h^2(w)}{M^2(w)}
=\frac{4w h^2(w)}{M(w)}\left(-\frac{d}{dw}\left\{\frac{1}{h(w)}\right\} 
+\frac{1}{4M(w)}\right).
\end{align*}
First we show that $h(w) \geq 0$ for all $w \geq 0$.
Suppose to the contrary that $ h(w)< 0 $ 
for some $ w_{0} $. Then  $ h(w) < 0 $ for all $ w \geq
w_{0} $ since $ h'(w) $ should be negative.
For all $ w > w_{0} $, the inequality
\begin{equation} \label{int-1}
\frac{d}{dw}\left( \frac{1}{h(w)}\right) \geq \frac{1}{4M(w)}
\end{equation}
should be satisfied.
Integrating both sides of \eqref{int-1}
from $ w_{0} $ to $ w^{*} $ leads to
\begin{equation*} 
\frac{1}{h(w^{*})}-\frac{1}{h(w_{0})} \geq
 \frac{1}{4}\int_{w_{0}}^{w^{*}}
M^{-1}(t)dt.
\end{equation*}
As $ w^{*} \rightarrow \infty $, the right-hand side of above inequality 
tends to infinity, and this provides a contradiction since the left-hand
 side is less than $ -1/h(w_{0}) $. Thus we have
$ g(w) \geq 0 $ for all $w$.

Similarly we can show that $ g(w) \leq 0 $ for all $w$. 
It follows that $ h(w) $ is zero
 for all $ w $, which implies that $ k(w) \equiv 0 $ for all $ w $.
This completes the proof.  
\end{proof}   
Combining Theorem \ref{perm} and \cite{Brown1971}'s sufficient condition
 for admissibility, we see that
 a quasi-admissible estimator of the form \eqref{pseudo}
is admissible
if it is
truly generalized Bayes,
that is,
there exists a nonnegative measure $\nu$ such that
\begin{align} \label{int-2}
 m(\|x\|^2)=\int_{R^p}\exp(-\|x-\theta\|^2/2)\nu(d\theta).
\end{align}
However, for given $m$, it is often quite difficult to determine
whether $m$ has such an integral form.
Furthermore, even if we find that $m$ does not have an integral form
like \eqref{int-2}, that is, 
the estimator is inadmissible,
it is generally very difficult to find an estimator dominating 
the inadmissible estimator.

The function $m$ for the James-Stein estimator is $m(w)=w^{2-p}$,
which satisfies the assumptions of Theorem \ref{perm}.
\begin{corollary}
The James-Stein estimator is quasi-admissible. 
\end{corollary}
In this case, it is not difficult to find 
an estimator dominating $\delta_{JS}$ by taking  positive part.
But it is not easy to find a large class of estimators
dominating the James-Stein estimator.
In Section 3, we introduce  an elegant sufficient condition for
domination over $\delta_{JS}$ proposed by
\cite{Kubokawa1994}.
%
%

\section{A class of quasi-admissible estimators improving upon the
 James-Stein estimator}

In this section, we introduce \cite{Kubokawa1994}'s sufficient condition
for improving upon the James-Stein estimator
and propose  a class of smooth quasi-admissible estimators satisfying it.

\cite{Kubokawa1994} showed that
if $ \lim_{w \rightarrow \infty}\phi(w)=p-2 $
  then the difference of risk functions between $\delta_{JS}$ and 
$\delta_\phi$ can be written as
\begin{align}
& R(\theta,\delta_{JS})-R(\theta,\delta_{\phi}) \nonumber \\
&   = 
 2\int_{0}^{\infty} \phi'(w)\left(  \phi(w)-(p-2) 
+ \frac{2f(w,\lambda)}{\int_{0}^{w}y^{-1} f_{p}(y,\lambda)dy} \right) 
\label{kubo-identity}\\
& \qquad \qquad \times \int_{0}^{w}y^{-1} f_{p}(y,\lambda)dy  dw , \nonumber
\end{align}
where $ \lambda=\| \theta \|^{2} $ and $ f_{p}(x;\lambda) $ denotes a
density of a non-central chi-square distribution with $ p $ degrees of
freedom and non-centrality parameter $ \lambda $. 

Moreover by the inequality
\[
 f_{p}(w;\lambda)/ \int_{0}^{w}y^{-1} f_{p}(y;\lambda)dy \geq 
 f_{p}(w)/ \int_{0}^{w}y^{-1} f_{p}(y)dy,
\]
where $ f_{p}(y)=f_{p}(y;0) $, which can be shown by the correlation
 inequality, 
we have
\begin{align*}
  & R(\theta,\delta_{JS})-R(\theta,\delta_{\phi}) \nonumber \\
  & \geq 2 \int_{0}^{\infty} \phi'(w)
 \left( \phi(w)-\phi_{0}(w)\right)
\left( \int_{0}^{w}y^{-1} f_{p}(y,\lambda)dy\right) dw, \label{saigo1}
\end{align*}
where 
\[
\phi_{0}(w)=p-2+2 f_{p}(w)/ \int_{0}^{w}y^{-1} f_{p}(y)dy. 
\]
Since  $ \phi_{0} = \phi_{K} $ given in \eqref{kubo-phi}
by an integration by parts and  $ \lim_{w \to \infty} 
\phi_K(w) = p-2$,
we have the following result.
\begin{thm}[\cite{Kubokawa1994}] \label{tatsuya}
If $ \phi(w) $ is nondecreasing and within $[\phi_K(w), p-2]$ for any 
$w \ge  0$, then
$ \delta_{\phi}(X) $ of form \eqref{estiortho} dominates the
James-Stein estimator.
\end{thm}
The assumption of the theorem above is satisfied by
$ \phi_K$ and $\phi_{JS}^+$. 
By \eqref{kubo-identity}, we see that 
the risk difference at $\|\theta\|=0$ 
between $\delta_{JS}$ and $ \delta_\phi$, the limit
of which is $p-2$, is given by
\begin{align} \label{risk-no-sa}
%
2 \int_{0}^{\infty} \phi'(w)
 \left( \phi(w)-\phi_{K}(w)\right)
\left( \int_{0}^{w}y^{-1} f_{p}(y)dy\right) dw.
\end{align}
Hence $\delta_K$ 
does not improve upon $ \delta_{JS}(X) $ at $ \| \theta \|=0 $,
although
$ \delta_{K}(X) $ is an admissible generalized Bayes estimator 
and thus smooth.
%
On the other hand, $ \delta_{JS}^+(X) $ 
improves on $ \delta_{JS}(X) $ 
significantly at the origin,
but it is not analytic and is thus inadmissible by \cite{Brown1971}'s
complete class theorem.
Therefore a more challenging  problem is to find admissible
estimators dominating the James-Stein estimator especially at 
$ \|\theta \|=0 $. 

In this paper, we 
propose a class of estimators
\begin{equation} \label{marusui}
 \delta_{\alpha}(X)=\left( 1-\phi_{\alpha}(\| X \|^{2})/\| X \|^{2} 
\right)X 
\end{equation}
where
\begin{align}
 \phi_{\alpha}(w) &= w\frac{\int_{0}^{1}\lambda^
{ \alpha(p/2-1)}\exp(-w \alpha \lambda/2
)d\lambda}{\int_{0}^{1}\lambda^
{ \alpha(p/2-1)-1}\exp(-w \alpha \lambda/2)d\lambda}   \nonumber  \\
&= p-2-\frac{2\exp(-w\alpha/2)}{\alpha\int_{0}^{1}\lambda^
{ \alpha(p/2-1)-1}\exp(-w \alpha \lambda/2)d\lambda} \nonumber \\
&= p-2-\frac{2}{\alpha\int_{0}^{1}(1-\lambda)^
{ \alpha(p/2-1)-1}\exp(w \alpha \lambda/2)d\lambda} .
\label{djs} 
\end{align}
The main theorem of this paper is as follows.
\begin{thm} \label{maru}
\begin{enumerate}
 \item $ \delta_{\alpha}(X)  $ dominates $ \delta_{JS}(X) $ 
for $ \alpha \geq 1 $.
\item The risk of $ \delta_{\alpha}(X)  $ 
for $\alpha>1$ at $\|\theta\|=0$ 
is strictly less than the risk of the James-Stein estimator at
       $\|\theta\|=0$. 
\item $ \delta_{\alpha}(X)  $ approaches 
the positive-part James-Stein estimator 
when $ \alpha $ tends to infinity, that is,
\[ \lim_{\alpha \rightarrow \infty}\delta_{\alpha}(X)
=\delta^+_{JS}(X). \]
\end{enumerate}
\end{thm}
Clearly  $ \delta_{1}(X)=\delta_{K}(X) $.
The class of $ \delta_\alpha $ with $ \alpha \ge 1$ is a bridge between
$\delta_K$ which is admissible and $ \delta_{JS}^+$ which is
inadmissible.
\begin{proof} \
[part 1]  We shall verify that $ \phi_{\alpha}(w) $ for $ \alpha \geq 1 $
satisfies assumptions in Theorem \ref{tatsuya}.
Applying the Taylor expansion to a part of \eqref{djs}, we have
\begin{align*}
& \frac{\alpha}{2}\int_{0}^{1}(1-\lambda)^
{ \alpha(p/2-1)-1}\exp(w \alpha \lambda/2)d\lambda \\
&=\sum_{i=0}^\infty w^i \prod_{j=0}^i (p-2+2j/\alpha)^{-1}
 =\psi(\alpha,w) \quad (say.)
\end{align*}
As $ \psi(\alpha,w) $ is increasing in $ w $, 
$ \phi_{\alpha}(w) $ is increasing in $ w $.
As $ \lim_{w \rightarrow \infty}\psi(\alpha,w)=\infty, $ for any 
$\alpha \geq 1$, 
it is clear that
$ \lim_{w \rightarrow \infty}\phi_{\alpha}(w)=p-2 $.
In order to show that $ \phi_{\alpha}(w) \geq \phi_{K}(w)
=\phi_{1}(w) $ for $ \alpha \geq 1 $, we have only to check that 
$ \psi(\alpha,w) $ is increasing in $ \alpha $.
It is easily verified because
the coefficient of each term of $\psi(\alpha,w)$ 
is increasing in $ \alpha $.
We have thus proved the theorem.

\medskip

[part 2]
Since $\phi_\alpha$ for $\alpha>1$ is strictly greater than
$\phi_K$ and strictly increasing in $w$,
the risk difference of $\delta_\alpha$ for $\alpha >1$ and $\delta_{JS}$
at $\|\theta\|=0$, which is given in
\eqref{risk-no-sa}, is strictly positive.

\medskip

[part 3] 
Since $\psi(\alpha,w)$ is increasing in $ \alpha $,
 it 
converges to 
\begin{align*}
(p-2)^{-1}\sum_{i=0}^{\infty}\left(\frac{w}{p-2} \right)^{i}
\end{align*}
by the monotone
convergence theorem when $\alpha$ goes to infinity. 
Considering two cases: $ w <( \geq ) p-2 $,
we obtain $ \lim_{\alpha \rightarrow \infty}\phi_{\alpha}(w)= w $ if
$ w < p-2 $; $ =p-2 $ otherwise.
This completes the proof.
\end{proof}
The estimator $\delta_\alpha$ is expressed as
$X+\nabla \log m_\alpha(\|X\|^2)$ where
\begin{align*}
 m_\alpha(w)= \left\{ \int_{0}^{1}\lambda^
{ \alpha(p/2-1)-1}\exp\left(-\frac{\alpha w}{2}
 \lambda \right) d\lambda \right\}^{1/\alpha}.
\end{align*}
Since $  m_\alpha(w) \sim w^{2-p}$ for sufficiently large 
$w$ by Tauberian's theorem and $ 0<m_\alpha(0)<\infty $,
we have the following result  by Theorem \ref{perm}.
\begin{corollary}
 $\delta_{\alpha}(X) $ is quasi-admissible for $p \geq 3$. 
\end{corollary}
Needless to say,  we are extremely interested in determining whether  
or not $ \delta_{\alpha}(X) $ for $ \alpha > 1 $ is admissible. 
Since $\delta_{\alpha}(X) $ with $ \alpha > 1 $ is quasi-admissible,
it is admissible if it is 
 generalized Bayes, that is, there exists a measure $ \nu $ which
satisfies
\begin{equation*} 
\int_{\mathit{R^{p}}} 
\exp (-\| \theta -x\|^{2}/2) \nu(d\theta) =  \left( \int_{0}^{1}\lambda^
{ \alpha(p/2-1)-1}\exp(-\alpha \| X \|^{2} \lambda/2 )d\lambda \right)^{1/\alpha}.
\end{equation*}  
I have no idea on the way to construct such a measure $\nu$ so far.
Even if we find that there is no $\nu$, which implies $ \delta_\alpha$
is inadmissible, it is very difficult to find an estimator dominating
$\delta_\alpha$ for $\alpha>1$.
\bibliographystyle{acmtrans-ims}
\bibliography{../maru-bib}

\begin{thebibliography}{}

\bibitem{Blyth1951}
{\sc Blyth, C.~R.} (1951).
\newblock On minimax statistical decision procedures and their admissibility.
\newblock {\em Ann. Math. Statist.\/}~{\em 22}, 22--42.
\MR{MR0039966 (12,622f)}

\bibitem{Bock1988}
{\sc Bock, M.~E.} (1988).
\newblock Shrinkage estimators: pseudo-{B}ayes rules for normal mean vectors.
\newblock In {\em Statistical decision theory and related topics, IV, Vol.\ 1
  (West Lafayette, Ind., 1986)}. Springer, New York, 281--297.
\MR{MR927108 (89b:62106)}

\bibitem{Brown1971}
{\sc Brown, L.~D.} (1971).
\newblock Admissible estimators, recurrent diffusions, and insoluble boundary
  value problems.
\newblock {\em Ann. Math. Statist.\/}~{\em 42}, 855--903.
\MR{MR0286209 (44 \#3423)}

\bibitem{Brown1988}
{\sc Brown, L.~D.} (1988).
\newblock The differential inequality of a statistical estimation problem.
\newblock In {\em Statistical decision theory and related topics, IV, Vol.\ 1
  (West Lafayette, Ind., 1986)}. Springer, New York, 299--324.
\MR{MR927109 (89d:62010)}

\bibitem{Brown-Zhang2006}
{\sc Brown, L.~D.} {\sc and} {\sc Zhang, L.~H.} (2006).
\newblock Estimators for gaussian models having a block-wise structure.
\newblock Mimeo,\\ (available at: {\tt
  http://www-stat.wharton.upenn.edu/~lbrown/teaching/Shrinkage/}).

\bibitem{Guo-Pal1992}
{\sc Guo, Y.~Y.} {\sc and} {\sc Pal, N.} (1992).
\newblock A sequence of improvements over the {J}ames-{S}tein estimator.
\newblock {\em J. Multivariate Anal.\/}~\textbf{42},~2, 302--317.
\MR{MR1183849 (93h:62009)}

\bibitem{James-Stein1961}
{\sc James, W.} {\sc and} {\sc Stein, C.} (1961).
\newblock Estimation with quadratic loss.
\newblock In {\em Proc. 4th Berkeley Sympos. Math. Statist. and Prob., Vol. I}.
  Univ. California Press, Berkeley, Calif., 361--379.
\MR{MR0133191 (24 \#A3025)}

\bibitem{Kubokawa1991}
{\sc Kubokawa, T.} (1991).
\newblock An approach to improving the {J}ames-{S}tein estimator.
\newblock {\em J. Multivariate Anal.\/}~\textbf{36},~1, 121--126.
\MR{MR1094273 (92b:62070)}

\bibitem{Kubokawa1994}
{\sc Kubokawa, T.} (1994).
\newblock A unified approach to improving equivariant estimators.
\newblock {\em Ann. Statist.\/}~\textbf{22},~1, 290--299.
\MR{MR1272084 (95h:62011)}

\bibitem{Kuriki-Takemura2000}
{\sc Kuriki, S.} {\sc and} {\sc Takemura, A.} (2000).
\newblock Shrinkage estimation towards a closed convex set with a smooth
  boundary.
\newblock {\em J. Multivariate Anal.\/}~\textbf{75},~1, 79--111.
\MR{MR1787403 (2001h:62034)}

\bibitem{Li-Kuo1982}
{\sc Li, T.~F.} {\sc and} {\sc Kuo, W.~H.} (1982).
\newblock Generalized {J}ames-{S}tein estimators.
\newblock {\em Comm. Statist. A---Theory Methods\/}~\textbf{11},~20,
  2249--2257.
\MR{MR678683 (84b:62078)}

\bibitem{Rukhin1995}
{\sc Rukhin, A.~L.} (1995).
\newblock Admissibility: {S}urvey of a concept in progress.
\newblock {\em International Statistical Review\/}~{\em 63}, 95--115.

\bibitem{Stein1956}
{\sc Stein, C.} (1956).
\newblock Inadmissibility of the usual estimator for the mean of a multivariate
  normal distribution.
\newblock In {\em Proceedings of the Third Berkeley Symposium on Mathematical
  Statistics and Probability, 1954--1955, vol. I}. University of California
  Press, Berkeley and Los Angeles, 197--206.
\MR{MR0084922 (18,948c)}

\bibitem{Stein1974}
{\sc Stein, C.} (1974).
\newblock Estimation of the mean of a multivariate normal distribution.
\newblock In {\em Proceedings of the Prague Symposium on Asymptotic Statistics
  (Charles Univ., Prague, 1973), Vol. II}. Charles Univ., Prague, 345--381.
\MR{MR0381062 (52 \#1959)}

\end{thebibliography}


\end{document}